\documentclass[12pt,a4paper]{amsart}
\usepackage[foot]{amsaddr}
\usepackage[english]{babel}
\usepackage{amsmath,amssymb,amsthm,latexsym,MnSymbol}
\usepackage[utf8]{inputenc}
\usepackage[colorlinks=true]{hyperref}
\usepackage[a4paper, portrait]{geometry}

\usepackage[dvipsnames]{color}
\definecolor{green}{rgb}{0.0, 0.5, 0.05}
\definecolor{teal}{rgb}{0.0, 0.5, 0.5}
\definecolor{purple}{rgb}{0.5, 0.0, 0.5}
\definecolor{fuchsia}{rgb}{1.0, 0.0, 1.0}

\newcommand{\Z}{\mathbf{Z}}
\newcommand{\Q}{\mathbf{Q}}

\newcommand{\R}{\mathbf{R}}
\newcommand{\C}{\mathbf{C}}
\newcommand{\h}{\mathcal{H}}

\DeclareMathOperator{\im}{Im}
\DeclareMathOperator{\re}{Re}
\DeclareMathOperator{\SL}{SL}
\newcommand{\abcd}[4]{\begin{pmatrix}#1&#2\\#3& #4\end{pmatrix}}
\newcommand{\sabcd}[4]{\left(\begin{smallmatrix}#1&#2\\#3& #4\end{smallmatrix}\right)}
\newcommand{\bolda}{{\boldsymbol a}}
\newcommand{\boldb}{{\boldsymbol b}}
\newcommand{\boldc}{{\boldsymbol c}}

\newcommand{\boldu}{{\boldsymbol u}}
\newcommand{\boldv}{{\boldsymbol v}}
\newcommand{\boldw}{{\boldsymbol w}}
\newcommand{\boldx}{{\boldsymbol x}}

\newtheorem{thm}{Theorem}
\newtheorem*{thm*}{Theorem}
\newtheorem{lem}{Lemma}

\newtheorem*{cor*}{Corollary}

\theoremstyle{definition}

\theoremstyle{remark}

\newtheorem*{remarks*}{Remarks}

\newtheorem*{question*}{Question}

\begin{document}

\author[F. Brunault]{François Brunault}

\date{\today}

\address{ÉNS Lyon, Unité de mathématiques pures et appliquées, 46 allée d'Italie, 69007 Lyon, France}

\email{francois.brunault@ens-lyon.fr}
\urladdr{http://perso.ens-lyon.fr/francois.brunault}

\begin{abstract}
Borisov and Gunnells have proved that certain linear combinations of products of Eisenstein series are Eisenstein series themselves, in analogy with the Manin relations for modular symbols. We devise a new method for determining and proving such relations, by differentiating with respect to the parameters of the Eisenstein series.
\end{abstract}

\title[On the Borisov--Gunnells relations for Eisenstein series]{On the Borisov--Gunnells relations for products of Eisenstein series}



\thanks{The author was supported by the research project ``Motivic homotopy, quadratic invariants and diagonal classes'' (ANR-21-CE40-0015) operated by the French National Research Agency (ANR)}

\maketitle

\section{Introduction} \label{sec:intro}

The Eisenstein series considered here are defined as follows. We denote by $\h$ the upper half-plane. Let $k \geqslant 0$ be an integer, and $\tau \in \h$. For a complex number $z$, written as $z = x_1 \tau + x_2$ with $x_1, x_2 \in \R$, we introduce the following Eisenstein-Kronecker series
\begin{equation} \label{def Kk}
\mathcal{K}^{(k)}_z(\tau, s) = \sideset{}{'}\sum_{m,n \in \Z} \frac{e^{2\pi i(mx_2-nx_1)}}{(m\tau+n)^k |m\tau+n|^{2s}} \qquad (s \in \C),
\end{equation}
where $\sideset{}{'}\sum$ indicates that the term $(m,n) = (0,0)$ is omitted. This type of series, featuring an additive character in the numerator, was first investigated by Kronecker, see \cite[Chap.~VII--VIII]{Wei76} and~\cite{CS16}. The case $k=0$, $s=1$ leads to the Kronecker limit formulas, which have significant applications in number theory \cite{Sie80, DIT18}. The case $k=s=1$ is related to Bloch's work on the elliptic dilogarithm and $L$-functions of CM elliptic curves \cite{Blo00, DO07}, and his results paved the way to Beilinson's conjectures on special values of $L$-functions \cite{Nek94}.

In this article, we are interested in the case $k \geqslant 1$, $s=0$, which classically gives rise to holomorphic Eisenstein series of weight $k$. More precisely, for $k \geqslant 1$, the series $\mathcal{K}^{(k)}_z(\tau, s)$ converges for $\re(s) > 1 - \frac{k}{2}$ and extends to a holomorphic function of $s \in \C$ \cite[Chap.~VIII, \S 12--13]{Wei76}. We then define
\begin{equation} \label{def EK}
E^{(k)}_z(\tau) = - \frac{(k-1)!}{(-2\pi i)^k} \, \mathcal{K}^{(k)}_z(\tau, 0).
\end{equation}
It follows from the definition \eqref{def Kk} that the function $z \mapsto E^{(k)}_z(\tau)$ is periodic with respect to the lattice $\Z \tau + \Z$, so that $E^{(k)}_z(\tau)$ is well-defined for $z \in \C/(\Z\tau+\Z)$.

Now let $N \geqslant 1$ be an integer. To any parameter $\bolda = (a_1, a_2) \in (\Z/N\Z)^2$, we associate the function on $\h$
\begin{equation} \label{def EkN}
E^{(k;N)}_\bolda(\tau) = E^{(k)}_{(\tilde{a}_1 \tau + \tilde{a}_2)/N}(\tau),
\end{equation}
for any choice of representatives $\tilde{a}_1, \tilde{a}_2$ of $a_1, a_2$ in $\Z$. It is known (Lemma \ref{lem EkN}) that $E^{(k;N)}_\bolda$ belongs to the space of modular forms $M_k(\Gamma(N))$ (see Section \ref{subsec:notations} for the definition of modular forms), except in the case when $k=2$ and $\bolda=(0,0)$ (the function $E^{(2;N)}_{(0,0)}$ is modular of weight $2$ for $\SL_2(\Z)$, but not holomorphic).

As is well-known, the graded algebra of modular forms for $\SL_2(\Z)$ is generated by the Eisenstein series of weight $4$ and $6$ \cite[Proposition 4]{Zag08}. This automatically produces linear relations between products of Eisenstein series of varying weight. In arbitrary level, Khuri-Makdisi has shown that for $k \geqslant 2$ and $N \geqslant 3$, every modular form in $M_k(\Gamma(N))$ is a polynomial in the Eisenstein series $E^{(1;N)}_\bolda$ \cite[Theorems 3.5 and 5.1]{Khu12}. Under certain conditions, pairwise products of similar Eisenstein series $\tilde{s}^{(k)}_{a/N}$ with $a \in \Z/N\Z$ (see \eqref{eq skaN} for the definition) span the space of modular forms for $\Gamma_1(N)$ \cite{BG03}. This has applications to computing Fourier expansions of modular forms at the cusps \cite{Coh19, DN18}. Products of Eisenstein series also occur in Ramanujan's famous differential system \cite[Proposition 15]{Zag08}. Note that the Rankin-Cohen brackets of modular forms \cite[Section~5.2]{Zag08} also involve pairwise products of (quasi)modular forms.

Our main result establishes linear dependence relations among pairwise products of the Eisenstein series $E^{(k;N)}_\bolda$, continuing work of Borisov and Gunnells \cite{BG01a, BG01b,BG03}, Pa\textcommabelow{s}ol \cite{Pas06}, Khuri-Makdisi and Raji \cite{KR17} and Zhang \cite{Zha20}. In order to state the result, we introduce the following notations. Let $k_1, k_2 \geqslant 0$ be integers. For $u,v \in \C/(\Z\tau+\Z)$, we define the symbol
\begin{equation*}
X^{k_1} Y^{k_2} [u, v]_\tau = E^{(k_1+1)}_u(\tau) E^{(k_2+1)}_v(\tau).
\end{equation*}
We extend this definition to a homogeneous polynomial $P(X,Y) = \sum_{k_1 + k_2 = \ell} c_{k_1, k_2} X^{k_1} Y^{k_2}$ of degree $\ell$ in $\C[X,Y]$ by setting
\begin{equation*}
P [u, v]_\tau = \sum_{k_1+k_2 = \ell} c_{k_1, k_2} X^{k_1} Y^{k_2} [u, v]_\tau.
\end{equation*}
There are corresponding symbols for discrete parameters in $(\Z/N\Z)^2$. For $\bolda, \boldb \in (\Z/N\Z)^2$, we define the function on $\h$
\begin{equation*}
X^{k_1} Y^{k_2} [\bolda, \boldb] = E^{(k_1+1;N)}_\bolda E^{(k_2+1;N)}_\boldb.
\end{equation*}
We define similarly $P [\bolda, \boldb]$ for a homogeneous polynomial $P$ of degree $\ell$ in $\C[X,Y]$. Note that if $\bolda$ and $\boldb$ are non-zero, then $P[\bolda, \boldb] \in M_{\ell+2}(\Gamma(N))$. The following theorem provides non-trivial relations between these modular forms.

\begin{thm} \label{main thm}
Let $\tau \in \h$. For any integers $k \geqslant 2$ and $k_1, k_2 \geqslant 0$ such that $k=k_1+k_2+2$, and any $u, v, w \in \C/(\Z \tau + \Z)$ such that $u+v+w = 0$ and $u, v, w \neq 0$, we have
\begin{equation} \label{eq thm 1}
\begin{split}
& X^{k_1} Y^{k_2} [u, v]_\tau + (-X-Y)^{k_1} X^{k_2} [v, w]_\tau + Y^{k_1} (-X-Y)^{k_2} [w, u]_\tau \\
& \qquad = \frac{(-1)^{k_2+1}}{k_2+1} E^{(k)}_u(\tau) + \frac{(-1)^{k_1+1}}{k_1+1} E^{(k)}_v(\tau) + (-1)^{k_1+k_2+1} \frac{k_1! k_2!}{(k_1+k_2+1)!} E^{(k)}_w(\tau).
\end{split}
\end{equation}
In particular, for any weight $k \geqslant 2$, any level $N \geqslant 1$ and any $\bolda, \boldb, \boldc \in (\Z/N\Z)^2$ such that $\bolda + \boldb + \boldc = 0$ and $\bolda, \boldb, \boldc \neq 0$, we have the identity in $M_k(\Gamma(N))$
\begin{equation} \label{eq thm 2}
\begin{split}
& X^{k_1} Y^{k_2} [\bolda, \boldb] + (-X-Y)^{k_1} X^{k_2} [\boldb, \boldc] + Y^{k_1} (-X-Y)^{k_2} [\boldc, \bolda] \\
& \qquad = \frac{(-1)^{k_2+1}}{k_2+1} E^{(k;N)}_\bolda + \frac{(-1)^{k_1+1}}{k_1+1} E^{(k;N)}_\boldb + (-1)^{k_1+k_2+1} \frac{k_1! k_2!}{(k_1+k_2+1)!} E^{(k;N)}_\boldc.
\end{split}
\end{equation}
\end{thm}

Although the relation \eqref{eq thm 2} involving modular forms might be considered as the main object of interest, our proof of Theorem~\ref{main thm} requires the use of Eisenstein-Kronecker series with complex parameters. The main input is a differential relation satisfied by the series $E^{(k)}_z(\tau)$, seen as a function of the elliptic variable $z$. This allows to prove \eqref{eq thm 1} by induction on the weight $k$.

Let us compare our theorem in more detail with previously known results. In the case $k=2$, the identity \eqref{eq thm 1} was proved by Zhang \cite[Corollary 1.4.9]{Zha20} using a different method. The identity \eqref{eq thm 1} for $k=3$ was derived in a different way in \cite{BZ23} and played a crucial role there to compute certain regulator integrals for modular curves.

For integers $k,N \geqslant 1$, Borisov and Gunnells also defined Eisenstein series $\tilde{s}^{(k)}_{a/N}$ with parameter $a \in \Z/N\Z$ for the group $\Gamma_1(N)$ \cite[Section 2.1]{BG03}. They are related to $E^{(k;N)}_\bolda$ via the formula
\begin{equation} \label{eq skaN}
\tilde{s}^{(k)}_{a/N}(\tau) = - N^{k-1} E^{(k;N)}_{(a,0)}(N\tau),
\end{equation}
as can be deduced from Lemma \ref{lem Fourier Ekz}. Borisov and Gunnells \cite[Theorem 6.2]{BG03} show the identity \eqref{eq thm 2} for them, but without an explicit right-hand side. Also, our Eisenstein series are more general. The proof of Borisov and Gunnells follows a quite different method and uses toric modular forms, which were introduced in \cite{BG01a}. They also mention in \cite[Remark 6.4]{BG03} that their theorem can be obtained by differentiating with respect to a real parameter.

In the pioneering work \cite{Eis47}, Eisenstein defined his eponymous series with complex parameters, and the idea of differentiating them with respect to the parameters is already present there. The relation with $E^{(k)}_z(\tau)$ is the following. In the notations of Weil \cite[Chap.~VIII, \S 12]{Wei76}, we have $\mathcal{K}^{(k)}_z(\tau, s) = K_k(0, z, k+s)$, where $K_k(x,x_0,s)$ is a double series associated to the lattice $\Z\tau+\Z$ and depending on two complex parameters~$x$ and $x_0$. Eisenstein \cite[p.~223]{Eis47} defines series $(k,x)$ for $k \geqslant 1$ and $x \in \C$, which coincide with the other specialisation $K_k(x, 0, k)$ for $k \geqslant 3$ (in the cases $k=1$ and $k=2$, he uses a different summation process, leading to slightly different functions). Using the relation $\frac{\partial}{\partial x} (k, x) = -k \, (k+1, x)$ and taking also products of such series, he obtains (in an equivalent form) the differential equation satisfied by the $\wp$ function, anticipating by nearly fifteen years the work of Weierstra\ss{} \cite[p.~225]{Eis47}. One difference with the series $E^{(k)}_z(\tau)$ is that the weight of $E^{(k)}_z(\tau)$ decreases upon differentiation (see Lemma~\ref{lem diff Ek}).

The Eisenstein series considered in \cite{Pas06, KR17} are the series $K_k(x,0,k)$. Pa\textcommabelow{s}ol \cite[Theorem~3.1]{Pas06} obtains the analogue of Theorem \ref{main thm} for these series, by making use of differentiation with respect to $x$. Khuri-Makdisi and Raji \cite[Theorem 2.8]{KR17} prove these relations modulo Eisenstein series using the Rankin-Selberg method. After restricting to discrete parameters in $\frac{1}{N}(\Z\tau+~\Z)$, the series $K_k(0, x, k)$ and $K_k(x, 0, k)$ are related by a discrete Fourier transform (see \cite[Section~6]{BZ23}). It would be interesting to recover Theorem \ref{main thm} from \cite{Pas06} using this technique.

Finally, we note that our method allows to determine the precise form of the identity \eqref{eq thm 1}, without knowing it in advance (see Section~\ref{sec:explicit}).

\subsection{Notations and definitions.} \label{subsec:notations}
For $z \in \C$, we write $e(z) = \exp(2\pi iz)$. We denote by $\h = \{\tau \in \C : \im(\tau) > 0\}$ the upper half-plane. For $\tau \in \h$, we write $q = e(\tau)$ and $q^\alpha = e(\alpha\tau)$ for any $\alpha \in \R$. The group $\SL_2(\R)$ acts on $\h$ by homographies, $\sabcd{a}{b}{c}{d} \cdot \tau = \frac{a\tau+b}{c\tau+d}$.

The principal congruence subgroup of $\SL_2(\Z)$ of level $N \geqslant 1$ is
\begin{equation*}
\Gamma(N) = \biggl\{ \abcd{a}{b}{c}{d} \in \SL_2(\Z) : \abcd{a}{b}{c}{d} \equiv \abcd{1}{0}{0}{1} \bmod{N} \biggr\}.
\end{equation*}
A congruence subgroup of $\SL_2(\Z)$ is a subgroup containing $\Gamma(N)$ for some $N$. The congruence subgroup $\Gamma_1(N)$ is defined by
\begin{equation*}
\Gamma_1(N) = \biggl\{ \abcd{a}{b}{c}{d} \in \SL_2(\Z) : \abcd{a}{b}{c}{d} \equiv \abcd{1}{*}{0}{1} \bmod{N} \biggr\}.
\end{equation*}

Given a congruence subgroup $\Gamma$ of $\SL_2(\Z)$, a modular form of weight $k \geqslant 1$ for $\Gamma$ is a holomorphic function $f : \h \to \C$ satisfying the following two properties:
\begin{enumerate}

\item[$\bullet$] (Modularity) For every $\gamma = \sabcd{a}{b}{c}{d} \in \Gamma$ and $\tau \in \h$, we have $f(\gamma \tau) = (c\tau+d)^k f(\tau)$.

\item[$\bullet$] (Holomorphy at infinity) For every $\gamma = \sabcd{a}{b}{c}{d} \in \SL_2(\Z)$, the function $\tau \mapsto (c\tau+d)^{-k} f(\gamma \tau)$ is bounded as $\im(\tau) \to +\infty$.

\end{enumerate}
We denote by $M_k(\Gamma)$ the complex vector space of modular forms of weight $k$ for $\Gamma$. There is a direct sum decomposition $M_k(\Gamma) = S_k(\Gamma) \oplus \mathcal{E}_k(\Gamma)$, where $S_k(\Gamma)$ is the space of cusp forms and $\mathcal{E}_k(\Gamma)$, called the Eisenstein subspace, is spanned by Eisenstein series \cite[Corollary 8.2.6]{CS17}.

For an integer $k \geqslant 2$, we denote by $V_{k-2}$ the vector space of homogeneous polynomials of degree $k-2$ in $\C[X,Y]$.

\bigskip

\textbf{Acknowledgements.} I would like to thank the referee for their comments that enriched the article and led to improvements in the exposition. I am also grateful to Pierre Charollois, Kamal Khuri-Makdisi and Wadim Zudilin for their valuable feedback.

\section{Eisenstein series with complex parameters} \label{sec:continuous}

Let $k \geqslant 1$. It follows from the definition \eqref{def Kk} that $\mathcal{K}^{(k)}_{-z}(\tau, s) = (-1)^k \mathcal{K}^{(k)}_z(\tau, s)$ for $\re(s) > 1 - \frac{k}{2}$. By analytic continuation, this is true for all $s \in \C$, and taking $s=0$, we get
\begin{equation} \label{eq parity}
E^{(k)}_{-z}(\tau) = (-1)^k E^{(k)}_z(\tau).
\end{equation}

We will make use of the Fourier expansion of $E^{(k)}_z(\tau)$ with respect to $\tau$ (see \cite[Lemma 34]{BZ23}).

\begin{lem} \label{lem Fourier Ekz}
Let $k \geqslant 1$ be an integer and $z \in \C$, written as $z = x_1 \tau + x_2$ with $x_1,x_2 \in \R$. Assume that $z \not\in \Z \tau + \Z$ in the case when $k=2$. Then
\begin{equation} \label{eq qexp Ekz}
E^{(k)}_z(\tau) = a_0(E^{(k)}_z) - \sum_{\substack{ \mu \geqslant 1 \\ \nu \in x_1+\Z \\ \nu > 0}} e(\mu x_2) \nu^{k-1} q^{\mu \nu} + (-1)^{k+1} \sum_{\substack{ \mu \geqslant 1 \\ \nu \in -x_1+\Z \\ \nu > 0}} e(-\mu x_2) \nu^{k-1} q^{\mu \nu},
\end{equation}
with
\begin{align*}
a_0(E^{(1)}_z) & = \begin{cases} 0 & \text{if } x_1, x_2 \in \Z, \\
-\frac12 \frac{1+e(x_2)}{1-e(x_2)} & \text{if } x_1 \in \Z \text{ and } x_2 \not\in \Z, \\
\{x_1\} - \frac12 & \text{if } x_1 \not\in \Z, \end{cases} \\
a_0(E^{(k)}_z) & = \frac{B_k(\{x_1\})}{k} \qquad (k \geqslant 2),
\end{align*}
where $B_k(t)$ is the $k$th Bernoulli polynomial and $\{ \,\cdot\, \}$ stands for the fractional part.
\end{lem}

We now discuss the Eisenstein-Kronecker series $E^{(k)}_z(\tau)$ when $z$ is an $N$-torsion point of the complex torus $\C/(\Z \tau+\Z)$. It will be convenient to view $E^{(k)}_z(\tau)$ as a function of two real variables by means of the decomposition $z = x_1 \tau + x_2$ with $x_1, x_2 \in \R$. Explicitly, we define, for a row vector $\boldx = (x_1, x_2) \in \R^2$,
\begin{equation} \label{def Ekx}
\tilde{E}^{(k)}_{\boldx}(\tau) = E^{(k)}_{x_1 \tau + x_2}(\tau) \qquad (\tau \in \h).
\end{equation}
(This is the series denoted by $E^{(k)}_\boldx$ in \cite[Definition 33]{BZ23}.) This object behaves nicely with respect to the action of $\SL_2(\Z)$ \cite[Lemma 35]{BZ23}. Namely, for $\gamma = \sabcd{a}{b}{c}{d} \in \SL_2(\Z)$, we have
\begin{equation} \label{eq modularity Ekx}
\tilde{E}^{(k)}_\boldx(\gamma \tau) = (c\tau+d)^k \tilde{E}^{(k)}_{\boldx \gamma}(\tau).
\end{equation}

\begin{lem} \label{lem EkN}
Let $\bolda = (a_1, a_2) \in (\Z/N\Z)^2$, with $\bolda \neq (0,0)$ in the case when $k=2$. Then $E^{(k;N)}_\bolda$ belongs to $M_k(\Gamma(N))$.
\end{lem}

\begin{proof}
Note that $\tilde{E}^{(k)}_\boldx$ depends only on the class of $\boldx$ in $(\R/\Z)^2$. Using the definition \eqref{def EkN} of $E^{(k;N)}_\bolda$ and the identity \eqref{eq modularity Ekx}, one checks that for $\gamma = \sabcd{a}{b}{c}{d} \in \SL_2(\Z)$,
\begin{equation*}
E^{(k;N)}_\bolda(\gamma \tau) = (c\tau+d)^k E^{(k;N)}_{\bolda \gamma}(\tau).
\end{equation*}
Since a row vector in $(\Z/N\Z)^2$ is invariant under right multiplication by $\Gamma(N)$, it follows that the function $E^{(k;N)}_\bolda$ is modular of weight $k$ for $\Gamma(N)$. Moreover, Lemma \ref{lem Fourier Ekz} shows that $E^{(k;N)}_\bolda$ is holomorphic at infinity.
\end{proof}

\section{A differential relation for Eisenstein series} \label{sec:differential}

In Sections \ref{sec:differential} and \ref{sec:explicit}, we will work with a fixed $\tau \in \h$ and investigate the Eisenstein series as a function of the elliptic variable $z \in \C$. To lighten the notations, we will write $E^{(k)}_z = E^{(k)}_z(\tau)$ and $P[u,v] = P[u,v]_\tau$.

One crucial ingredient of the proof of Theorem \ref{main thm} is the following differential property of $E^{(k)}_z$ as a function of $z$.

\begin{lem} \label{lem diff Ek}
For any $k \geqslant 1$, the function $z \mapsto E^{(k)}_z$ is smooth on $\C \setminus (\Z \tau + \Z)$.
On this domain, we have
\begin{equation} \label{eq diff Ek}
-(\tau-\overline{\tau}) \frac{\partial}{\partial \overline{z}} E^{(1)}_z = 1, \qquad -(\tau-\overline{\tau}) \frac{\partial}{\partial \overline{z}} E^{(k)}_z = (k-1) E^{(k-1)}_z \qquad (k \geqslant 2).
\end{equation}
\end{lem}

\begin{proof}
The function $z \mapsto E^{(k)}_z$ is smooth away from $\Z \tau + \Z$ by \cite[Lemma 38]{BZ23}. Assume first that $k \geqslant 2$. The Eisenstein-Kronecker series $\mathcal{K}^{(k)}_z(\tau, s)$ satisfies the following differential property (compare with \cite[last equation of Lemma~1.4]{BK10}). For $z = x_1 \tau + x_2$ with $x_1, x_2 \in \R$, and $\lambda = m \tau + n$, we have
\begin{equation*}
e(mx_2 - nx_1) = e \Bigl( \frac{\lambda \bar{z} - \bar{\lambda} z}{\tau - \bar{\tau}} \Bigr).
\end{equation*}
A direct computation shows that
\begin{equation*}
\frac{\partial}{\partial \overline{z}} \mathcal{K}^{(k)}_z(\tau, s) = \frac{\pi}{\im(\tau)} \mathcal{K}^{(k-1)}_z(\tau, s).
\end{equation*}
This holds for $\re(s) > \frac{3-k}{2}$, and by analytic continuation for all $s \in \C$. Taking $s=0$ and using the definition \eqref{def EK} of $E^{(k)}_z$, this implies the second formula of \eqref{eq diff Ek}.

For $k=1$, we will argue by using the Fourier expansion of $E^{(1)}_z$ given in Lemma \ref{lem Fourier Ekz}. Since $z \mapsto E^{(1)}_z$ is smooth on $\C \setminus (\Z \tau + \Z)$, the function $z \mapsto \frac{\partial}{\partial \overline{z}} E^{(1)}_z$ is continuous on this domain, and it suffices to prove the identity on the dense open subset $\{z \in \C : z = x_1 \tau + x_2, \, x_1 \not\in \Z\}$. Moreover $z \mapsto E^{(1)}_z$ is $(\Z \tau + \Z)$-periodic, so we may assume $0<x_1<1$. In this case
\begin{equation*}
E^{(1)}_z = x_1 - \frac12 - \sum_{\substack{m \geqslant 1 \\ n \geqslant 0}} e(mz) q^{mn} + \sum_{\substack{m \geqslant 1 \\ n \geqslant 1}} e(-mz) q^{mn}.
\end{equation*}
It follows that
\begin{equation*}
\frac{\partial}{\partial \overline{z}} E^{(1)}_z = \frac{\partial x_1}{\partial \overline{z}} = \frac{\partial}{\partial \overline{z}} \Bigl(\frac{z-\overline{z}}{\tau-\overline{\tau}}\Bigr) = -\frac{1}{\tau-\overline{\tau}}. \qedhere
\end{equation*}
\end{proof}

We now study the derivatives of the symbols $P[u, v]$. Let $k \geqslant 2$ be an integer. Recall that $V_{k-2}$ is the complex vector space of homogeneous polynomials $P(X,Y)$ of degree $k-2$.

\begin{lem} \label{lem diff Pab}
For every $P \in V_{k-2}$ and every $u, v \in \C \setminus (\Z \tau + \Z)$, we have
\begin{align*}
-(\tau-\overline{\tau}) \frac{\partial}{\partial \overline{u}} P[u, v] & = \frac{\partial P}{\partial X} [u, v] + P(0,1) E^{(k-1)}_v, \\
-(\tau-\overline{\tau}) \frac{\partial}{\partial \overline{v}} P[u, v] & = \frac{\partial P}{\partial Y} [u, v] + P(1,0) E^{(k-1)}_u.
\end{align*}
\end{lem}

\begin{proof}
The first formula to be proved is linear in $P$, so it suffices to consider $P = X^{k_1} Y^{k_2}$ with $k_1+k_2 = k-2$. In the case $k_1 \geqslant 1$, we have by Lemma \ref{lem diff Ek}
\begin{equation*}
-(\tau-\overline{\tau}) \frac{\partial}{\partial \overline{u}} E^{(k_1+1)}_u E^{(k_2+1)}_v = k_1 E^{(k_1)}_u E^{(k_2+1)}_v = k_1 X^{k_1-1} Y^{k_2} [u, v] = \frac{\partial P}{\partial X} [u, v].
\end{equation*}
This concludes the argument, since in this case $P(0,1)=0$. For $P = Y^{k-2}$, we have
\begin{equation*}
-(\tau-\overline{\tau}) \frac{\partial}{\partial \overline{u}} E^{(1)}_u E^{(k-1)}_v = E^{(k-1)}_v,
\end{equation*}
which is what we want since $\partial P/\partial X = 0$ and $P(0,1)=1$. The formula for the derivative of $P[u, v]$ with respect to $\overline{v}$ is proved similarly.
\end{proof}

\section{Determining the explicit form of the relations} \label{sec:explicit}

Let $k \geqslant 2$ be an integer. Consider the polynomial $P_{k_1, k_2} = X^{k_1} Y^{k_2}$ with $k_1 + k_2 = k-2$. We postulate the following shape of the Borisov--Gunnells 3-term relations:
\begin{equation} \label{eq Ansatz}
P_{k_1, k_2} [u, v] + Q_{k_1, k_2} [v, -u-v] + R_{k_1, k_2} [-u-v, u] = \alpha_{k_1, k_2} E^{(k)}_u + \beta_{k_1, k_2} E^{(k)}_v + \gamma_{k_1, k_2} E^{(k)}_{-u-v},
\end{equation}
Here $Q_{k_1, k_2}$ and $R_{k_1, k_2}$ are some unknown polynomials in $V_{k-2}$, and $\alpha_{k_1, k_2}$, $\beta_{k_1, k_2}$, $\gamma_{k_1, k_2}$ are some constants.

Differentiating \eqref{eq Ansatz} with respect to $\overline{u}$ and using Lemma \ref{lem diff Pab}, we have
\begin{equation} \label{eq du Ansatz}
\begin{split}
& \frac{\partial P_{k_1, k_2}}{\partial X} [u, v] + P_{k_1, k_2}(0,1) E^{(k-1)}_v - \frac{\partial Q_{k_1, k_2}}{\partial Y} [v, -u-v] - Q_{k_1, k_2}(1,0) E^{(k-1)}_v \\
& \quad - \frac{\partial R_{k_1, k_2}}{\partial X} [-u-v, u] - R_{k_1, k_2}(0,1) E^{(k-1)}_u + \frac{\partial R_{k_1, k_2}}{\partial Y} [-u-v, u] + R_{k_1, k_2}(1,0) E^{(k-1)}_{-u-v} \\
= \; & (k-1) \alpha_{k_1, k_2} E^{(k-1)}_u - (k-1) \gamma_{k_1, k_2} E^{(k-1)}_{-u-v}.
\end{split}
\end{equation}
Note that $\partial P_{k_1, k_2}/\partial X = k_1 X^{k_1-1} Y^{k_2} = k_1 P_{k_1-1,k_2}$. Here we use the convention that $P_{k_1, k_2}$, $Q_{k_1, k_2}$, $R_{k_1, k_2}$, $\alpha_{k_1, k_2}$, $\beta_{k_1, k_2}$ and $\gamma_{k_1, k_2}$ are zero whenever one of the indices $k_1, k_2$ is equal to $-1$.

In order for the identity \eqref{eq du Ansatz} to match the one corresponding to the indices $(k_1-1, k_2)$ in \eqref{eq Ansatz}, we should have
\begin{equation} \label{eq du conditions}
\begin{split}
-\frac{\partial Q_{k_1, k_2}}{\partial Y} & = k_1 Q_{k_1-1, k_2} \qquad \Bigl(\frac{\partial}{\partial Y}- \frac{\partial}{\partial X}\Bigr) R_{k_1, k_2} = k_1 R_{k_1-1, k_2} \\
(k-1) \alpha_{k_1, k_2} + R_{k_1, k_2}(0,1) & = k_1 \alpha_{k_1-1, k_2} \qquad Q_{k_1, k_2}(1,0) - P_{k_1, k_2}(0,1) = k_1 \beta_{k_1-1, k_2} \\
-(k-1) \gamma_{k_1, k_2} - R_{k_1, k_2}(1,0) & = k_1 \gamma_{k_1-1, k_2}.
\end{split}
\end{equation}
Similarly, differentiating \eqref{eq Ansatz} with respect to $\overline{v}$, we should have
\begin{equation} \label{eq dv conditions}
\begin{split}
\Bigl(\frac{\partial}{\partial X} - \frac{\partial}{\partial Y}\Bigr) Q_{k_1, k_2} & = k_2 Q_{k_1, k_2-1} \qquad - \frac{\partial R_{k_1, k_2}}{\partial X} = k_2 R_{k_1, k_2-1} \\
R_{k_1, k_2}(0,1) - P_{k_1, k_2}(1,0) & = k_2 \alpha_{k_1, k_2-1} \qquad (k-1) \beta_{k_1, k_2} + Q_{k_1, k_2}(1,0) = k_2 \beta_{k_1, k_2-1} \\
-(k-1) \gamma_{k_1, k_2} - Q_{k_1, k_2}(0,1) & = k_2 \gamma_{k_1, k_2-1}.
\end{split}
\end{equation}
We put $Q_{0,0} = R_{0,0} = 1$ and $\alpha_{0,0} = \beta_{0,0} = \gamma_{0,0} = -1$. For these values, the identities \eqref{eq du conditions} and \eqref{eq dv conditions} are then all satisfied for $k_1 = k_2 = 0$. Furthermore, note that for $k > 2$, a polynomial in $V_{k-2}$ is uniquely determined by its two partial derivatives. By induction, the polynomials $Q_{k_1, k_2}$ and $R_{k_1, k_2}$ with $(k_1, k_2) \neq (0,0)$ are therefore uniquely determined by the conditions \eqref{eq du conditions} and \eqref{eq dv conditions}. We find that the following polynomials fulfill these conditions:
\begin{equation*}
P_{k_1, k_2} = X^{k_1} Y^{k_2}, \qquad Q_{k_1, k_2} = (-X-Y)^{k_1} X^{k_2}, \qquad R_{k_1, k_2} = Y^{k_1} (-X-Y)^{k_2}.
\end{equation*}
Likewise, the numbers $\alpha_{k_1, k_2}$, $\beta_{k_1, k_2}$ and $\gamma_{k_1, k_2}$ with $(k_1, k_2) \neq (0,0)$ are uniquely determined by \eqref{eq du conditions} and \eqref{eq dv conditions}, and we find the following solutions:
\begin{equation*}
\alpha_{k_1, k_2} = \frac{(-1)^{k_2+1}}{k_2+1}, \qquad \beta_{k_1, k_2} = \frac{(-1)^{k_1+1}}{k_1+1}, \qquad \gamma_{k_1, k_2} = (-1)^{k_1+k_2+1} \frac{k_1! k_2!}{(k_1+k_2+1)!}.
\end{equation*}
We now have the precise form of the identity to be proved (equation \eqref{eq thm 1} in Theorem \ref{main thm}).

\section{Proof of the main theorem} \label{sec:proof}

We prove the identity \eqref{eq thm 1} of Theorem \ref{main thm} by induction on the weight $k$, using the differential properties of the series $E^{(k)}_z$ described in Section \ref{sec:differential}.

For $\tau \in \h$, let us denote by $H(u, v, \tau)$ the difference of the left-hand side and the right-hand side of \eqref{eq thm 1}, where $w$ is uniquely determined by the condition $u+v+w=0$. We are going to show that $H(u, v, \tau)$ does not depend on $u$ and $v$, and that the resulting function of $\tau$ is a holomorphic modular form of weight $k$ for $\SL_2(\Z)$. Finally, we will prove that this modular form is zero.

Let $k \geqslant 2$ be an integer such that Theorem \ref{main thm} holds for all weights less than $k$ (if $k=2$, this condition is empty). Let $k_1, k_2 \geqslant 0$ be such that $k=k_1+k_2+2$, and fix $v \in \C \setminus (\Z \tau + \Z)$. Note that $H(u, v, \tau)$ is $(\Z \tau + \Z)$-periodic with respect to $u$. The computation in Section \ref{sec:explicit} and the induction hypothesis show that the function $u \mapsto H(u, v, \tau)$ is holomorphic on the domain $u \not\equiv 0, -v \bmod{\Z \tau + \Z}$. We are going to show that this function is bounded. Since a bounded elliptic function is constant, this will imply that $H(u, v, \tau)$ does not depend on $u$.

To this end, we need to control the behaviour of $E^{(k)}_z$ as the elliptic variable $z \in \C$ tends to a point of the lattice $\Z \tau + \Z$. By periodicity, it suffices to look at what happens as $z \to 0$.

\begin{lem} \label{lem Ek origin}
$\bullet$ If $k \geqslant 3$, then $z \mapsto E^{(k)}_z$ is continuous at the points of $\Z \tau + \Z$.

$\bullet$ If $k=2$, then for $z \to 0$, $z \neq 0$, with $z=x_1\tau+x_2$, $x_1, x_2 \in \R$, we have
\begin{equation*}
E^{(2)}_z = -2 G_2(\tau) + \frac{x_1}{2\pi i z} + \varepsilon_2(z),\end{equation*}
where $\displaystyle G_2(\tau) = -\frac1{24} + \sum_{m,n \geqslant 1} n q^{mn}$, and $\varepsilon_2(z)/z$ is bounded as $z \to 0$. In particular $E^{(2)}_z(\tau)$ is bounded as $z \to 0$ (but is not continuous at $z=0$).

$\bullet$ If $k=1$, then for $z \to 0$, $z \neq 0$, we have
\begin{equation*}
E^{(1)}_z = \frac{1}{2\pi iz} + \varepsilon_1(z),
\end{equation*}
where $\varepsilon_1(z)/z$ is bounded as $z \to 0$.
\end{lem}

\begin{proof}
$\bullet$ If $k \geqslant 3$, then the series defining $\mathcal{K}^{(k)}_z(\tau, 0)$ is normally convergent, hence defines a continuous function of $z$ on $\C$.

$\bullet$ Assume $k=2$. Since $E^{(2)}_{-z} = E^{(2)}_z$ by \eqref{eq parity}, it suffices to prove the asymptotics as $z = x_1 \tau + x_2$ tends to $0$ with $x_1 \geqslant 0$. In what follows, we assume $x_1 \in [0,\frac12]$. By Lemma \ref{lem Fourier Ekz}, we have
\begin{equation} \label{eq E2}
E^{(2)}_z = A(x_1, x_2) - x_1 \sum_{m \geqslant 1} e(mx_2) q^{mx_1},
\end{equation}
where
\begin{equation*}
A(x_1, x_2) = \frac{B_2(x_1)}{2} - \sum_{m,n \geqslant 1} e(mx_2) (n+x_1) q^{m(n+x_1)} - \sum_{m,n \geqslant 1} e(-mx_2) (n-x_1) q^{m(n-x_1)}
\end{equation*}
is a $C^\infty$ function of $x_1, x_2$, whose value at $x_1 = x_2 = 0$ is $-2 G_2(\tau)$. The second term of \eqref{eq E2}, to be considered only when $x_1>0$, is
\begin{equation*}
-x_1 \frac{e(z)}{1-e(z)} = \frac{x_1}{2\pi i z} + \varepsilon(z),
\end{equation*}
where $\varepsilon(z)/z$ is bounded as $z \to 0$. This gives the desired estimate.

$\bullet$ Assume $k=1$. Since $E^{(1)}_{-z} = -E^{(1)}_z$ by \eqref{eq parity}, it suffices to prove the asymptotics as $z = x_1 \tau + x_2$ tends to $0$ with $x_1 \geqslant 0$. We again assume $x_1 \in [0,\frac12]$. We use the Fourier expansion of $E^{(1)}_z$ from Lemma \ref{lem Fourier Ekz}, distinguishing the cases $x_1>0$ and $x_1=0$. If $x_1>0$, we have
\begin{align*}
E^{(1)}_z & = x_1 - \frac12 - \sum_{\substack{m \geqslant 1 \\ n \geqslant 0}} e(mx_2) q^{m(n+x_1)} + \sum_{\substack{m \geqslant 1 \\ n \geqslant 1}} e(-mx_2) q^{m(n-x_1)}.
\end{align*}
The terms of the series with $m, n \geqslant 1$ define a $C^\infty$ function of $x_1, x_2$, whose value at $x_1 = x_2 = 0$ is zero. The remaining terms in $E^{(1)}_z$ are
\begin{equation} \label{eq E1 origin}
x_1 - \frac12 - \sum_{m \geqslant 1} e(mz) = x_1 - \frac12 - \frac{e(z)}{1-e(z)} = \frac{1}{2\pi i z} + \varepsilon(z),
\end{equation}
where $\varepsilon(z)/z$ is bounded as $z \to 0$. In the case $x_1 = 0$, we have $z = x_2$, so that
\begin{align*}
E^{(1)}_z & = -\frac12 \cdot \frac{1+e(z)}{1-e(z)} - \sum_{m, n \geqslant 1} e(mz) q^{mn} + \sum_{m,n \geqslant 1} e(-mz) q^{mn}.
\end{align*}
Again, the two series over $m,n \geqslant 1$ define a $C^\infty$ function of $z$ whose value at $z = 0$ is zero. The remaining term is
\begin{equation*}
-\frac12 \cdot \frac{1+e(z)}{1-e(z)},
\end{equation*}
which coincides with \eqref{eq E1 origin} since $x_1=0$, and thus has the required shape.
\end{proof}

Let us go back to $H(u, v, \tau)$ and its behaviour as $u$ tends to $0$. Recall that $v \in \C \setminus (\Z \tau + \Z)$ is fixed. By Lemma \ref{lem diff Ek}, the function $u \mapsto E^{(\ell)}_{-u-v}$ is bounded as $u \to 0$ for any $\ell \geqslant 1$. So is the function $u \mapsto E^{(\ell)}_u$ for any $\ell \geqslant 2$, by Lemma \ref{lem Ek origin}. Therefore, it suffices to look at the contributions of $E^{(1)}_u$ in the expression of $H(u, v, \tau)$. This series appears only when $k_1 = 0$, and the contribution in this case is
\begin{equation} \label{eq u=0}
\begin{split}
Y^{k-2} [u, v] + (-X)^{k-2} [-u-v, u] & = E^{(1)}_u E^{(k-1)}_v + (-1)^{k-2} E^{(k-1)}_{-u-v} E^{(1)}_u \\
& = E^{(1)}_u (E^{(k-1)}_v - E^{(k-1)}_{u+v}),
\end{split}
\end{equation}
where we use the identity $E^{(k-1)}_{-z} = (-1)^{k-1} E^{(k-1)}_z$ from \eqref{eq parity}. By Lemma \ref{lem Ek origin}, the quantity $u E^{(1)}_u$ is bounded as $u \to 0$. Moreover, the function $z \mapsto E^{(k-1)}_z$ is smooth at $z = v$ by Lemma \ref{lem diff Ek}. It follows that the expression \eqref{eq u=0} is bounded as $u \to 0$.

The same analysis can be carried out when $u$ tends to $-v$. The only contributing terms in $H(u, v, \tau)$ are those involving $E^{(1)}_{-u-v}$, namely
\begin{equation*}
\begin{split}
(-X)^{k_1} X^{k_2} [v, -u-v] + Y^{k_1} (-Y)^{k_2} [-u-v, u] & = (-1)^{k_1} E^{(k-1)}_v E^{(1)}_{-u-v} + (-1)^{k_2} E^{(1)}_{-u-v} E^{(k-1)}_u \\
& = (-1)^{k_2} E^{(1)}_{-u-v} (E^{(k-1)}_u - E^{(k-1)}_{-v}),
\end{split}
\end{equation*}
which is again bounded as $u \to -v$. Since a bounded elliptic function is constant, we deduce that $H(u, v, \tau)$ does not depend on~$u$. Because of the symmetry $P(X, Y) [u, v] = P(Y, X) [v, u]$, we have $H(u, v, \tau) = H(v, u, \tau)$, so that $H(u, v, \tau)$ does not depend on $v$ either.

At this point, we need to investigate the dependence on the modular variable $\tau$. This requires considering additional symbols depending on real variables as in \eqref{def Ekx}. For a polynomial $P \in V_{k-2}$ and two real vectors $\boldu = (u_1, u_2)$ and $\boldv = (v_1, v_2)$, we set
\begin{equation*}
P[\boldu, \boldv](\tau) = P[u_1 \tau + u_2, v_1 \tau + v_2]_\tau \qquad (\tau \in \h).
\end{equation*}
For $\gamma = \sabcd{a}{b}{c}{d} \in \SL_2(\Z)$, we have thanks to \eqref{eq modularity Ekx} the transformation formula
\begin{equation} \label{eq modularity Puv}
P[\boldu, \boldv](\gamma \tau) = (c\tau+d)^k P[\boldu \gamma, \boldv \gamma](\tau).
\end{equation}

Note that the function $f(\tau) = H(u_1 \tau + u_2, v_1 \tau + v_2, \tau)$ does not depend on the choice of vectors $\boldu = (u_1, u_2)$ and $\boldv = (v_1, v_2)$ in $\R^2$. It follows from \eqref{eq modularity Puv} that for any $\gamma = \sabcd{a}{b}{c}{d} \in \SL_2(\Z)$, we have
\begin{equation*}
f(\gamma \tau) = (c\tau+d)^k f(\tau).
\end{equation*}
Moreover, Lemma \ref{lem Fourier Ekz} shows that $f$ is holomorphic at infinity. Thus $f$ belongs to $M_k(\SL_2(\Z))$.

Finally, we specialise the parameters $u$, $v$ to points of infinite order in $\C/(\Z \tau + \Z)$. In the Fourier expansion \eqref{eq qexp Ekz} of $E^{(k)}_z$, if we except the constant term, then all the exponents of $q$ are of the form $m x_1 + n$ with $m \in \Z \setminus \{0\}$ and $n \in \Z$. Now fix two vectors $\boldu, \boldv$ in $\R^2$ such that $0<u_1, v_1<1$ and $(1, u_1, v_1)$ are linearly independent over $\Q$, and consider the $q$-expansions of both sides of \eqref{eq thm 1}. If we except the constant term, then all the exponents of $q$ in $X^{k_1} Y^{k_2} [\boldu, \boldv]$ are of the form $m u_1 + m' v_1 + n$ with $(m, m') \in \Z^2$, $(m,m') \neq (0,0)$, and thus cannot be integers. Let us show that the same is true for $f$. If we write $\boldw = -\boldu-\boldv$, then the exponents of $q$ coming from the second term in the left-hand side of \eqref{eq thm 1} are of the form $m v_1 + m' w_1 + n$ with $(m, m') \neq (0,0)$. Since $w_1 = -u_1-v_1$, these exponents have the required shape. The same holds for the third term, and for the right-hand side of \eqref{eq thm 1}, which proves the claim. But the Fourier expansion of $f$ involves only integer powers of $q$. This implies that $f$ is constant, and in fact $f = 0$ since $k \geqslant 2$. This finishes the proof of Theorem \ref{main thm}.

\section{Some open questions}

To conclude, we highlight several open questions suggested by our results. The assumptions in Theorem \ref{main thm} exclude the case $N=1$, and therefore do not cover Eisenstein series for $\SL_2(\Z)$. Is there an analogue of Theorem \ref{main thm} in this case? A natural idea would be to let $u,v \to 0$ in \eqref{eq thm 1}, but by Lemma \ref{lem Ek origin}, the function $E^{(1)}_z$ has a singularity at $z=0$, and $E^{(2)}_z$ is also discontinuous there.

Because of these assumptions, we also do not deal with the non-holomorphic Eisenstein series $E^{(2;N)}_{(0,0)}$ in this article. It would be interesting to extend Theorem \ref{main thm} to this case, perhaps in the context of quasimodular forms. In the same vein, one could investigate the non-holomorphic Eisenstein-Kronecker series $\mathcal{K}^{(k)}_z(\tau, s)$, obtained by taking integers $s \geqslant 1$ in the definition \eqref{def Kk}.

Another possible direction would be to find an analogue of Theorem \ref{main thm} for Rankin-Cohen brackets of Eisenstein series.

Finally, it would be interesting to explore the generalisation of Theorem \ref{main thm} to products of more than two Eisenstein series, as suggested in \cite[Remark 7.15]{BG03}.

\end{document}